\documentclass[11pt,a4paper]{article}

\usepackage[left=2.45cm, top=2.45cm,bottom=2.45cm,right=2.45cm]{geometry}
\usepackage{mathtools,amssymb,amsthm,mathrsfs,calc,graphicx,xcolor,cleveref,dsfont,tikz,pgfplots,bm, extarrows}
\usepackage[british]{babel}
\usepackage{amsfonts}              
\usepackage[T1]{fontenc}
\usetikzlibrary{positioning,arrows}
\numberwithin{equation}{section}

\newtheorem {theorem}{Theorem}

{\theoremstyle{definition}
	
}
{\theoremstyle{theorem}
	\newtheorem {remark}[theorem]{Remark}
	
}

\newcommand{\Cov}{\operatorname{Cov}}
\newcommand{\var}{\operatorname{var}}

\def\EE{\mathbb{E}}

\def\NN{\mathbb{N}}

\def\PP{\mathbb{P}}

\def\RR{\mathbb{R}}

\def\cC{\mathcal{C}}

\makeatletter
\let\@fnsymbol\@alph
\makeatother

\begin{document}
	
	\title{\bfseries {Approaching the coupon collector's problem\\ with group drawings via Stein's method}}
	
	\author{Carina Betken\footnotemark[1]\; and Christoph Th\"ale\footnotemark[2]}
	
	\date{}
	\renewcommand{\thefootnote}{\fnsymbol{footnote}}
	\footnotetext[1]{Ruhr University Bochum, Germany. Email: carina.betken@rub.de}
	\footnotetext[2]{Ruhr University Bochum, Germany. Email: christoph.thaele@rub.de}
	
	\maketitle
	
	\begin{abstract}
		\noindent  In this paper the coupon collector's problem with group drawings is studied. Assume there are $ n $ different coupons. At each time precisely $ s $ of the $ n $ coupons are drawn, where all choices are supposed to have equal probability. The focus lies on the fluctuations, as $n\to\infty$, of the number $Z_{n,s}(k_n)$ of coupons that have not been drawn in the first $k_n$ drawings. Using a size-biased coupling construction together with Stein's method for normal approximation, a quantitative central limit theorem for $Z_{n,s}(k_n)$ is shown for the case that $k_n={n\over s}(\alpha\log(n)+x)$, where $0<\alpha<1$ and $x\in\mathbb{R}$. The same coupling construction is used to retrieve a quantitative Poisson limit theorem in the boundary case $\alpha=1$, again using Stein's method. 
	\end{abstract}

\section{Introduction }
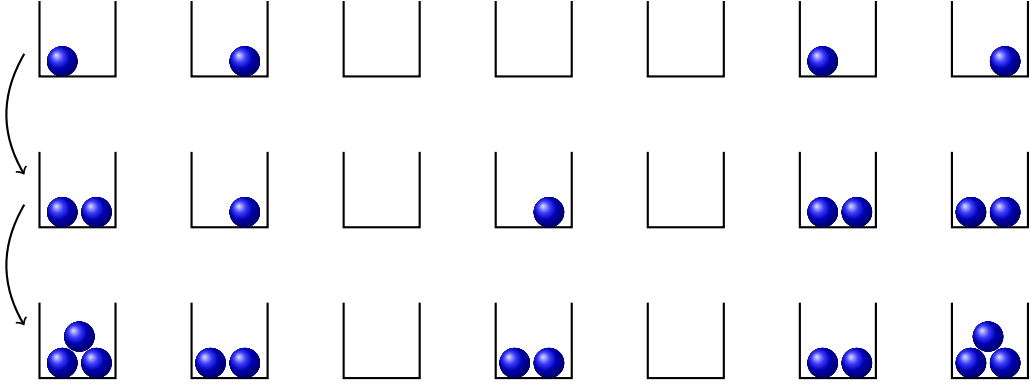
\begin{figure}[t]
 	\centering
\begin{tikzpicture}
	\draw[thick] (0,1) -- (0,0) -- (1,0) -- (1,1);
	\draw[thick] (2,1) -- (2,0) -- (3,0) -- (3,1);
	\draw[thick] (2+2,1) -- (2+2,0) -- (3+2,0) -- (3+2,1);
	\draw[thick] (2+4,1) -- (2+4,0) -- (3+4,0) -- (3+4,1);
	\draw[thick] (2+6,1) -- (2+6,0) -- (3+6,0) -- (3+6,1);
	\draw[thick] (2+8,1) -- (2+8,0) -- (3+8,0) -- (3+8,1);
	\draw[thick] (2+10,1) -- (2+10,0) -- (3+10,0) -- (3+10,1);
	
	\draw[thick] (0,1+2) -- (0,0+2) -- (1,0+2) -- (1,1+2);
	\draw[thick] (2,1+2) -- (2,0+2) -- (3,0+2) -- (3,1+2);
	\draw[thick] (2+2,1+2) -- (2+2,0+2) -- (3+2,0+2) -- (3+2,1+2);
	\draw[thick] (2+4,1+2) -- (2+4,0+2) -- (3+4,0+2) -- (3+4,1+2);
	\draw[thick] (2+6,1+2) -- (2+6,0+2) -- (3+6,0+2) -- (3+6,1+2);
	\draw[thick] (2+8,1+2) -- (2+8,0+2) -- (3+8,0+2) -- (3+8,1+2);
	\draw[thick] (2+10,1+2) -- (2+10,0+2) -- (3+10,0+2) -- (3+10,1+2);
	
	\draw[thick] (0,1+4) -- (0,0+4) -- (1,0+4) -- (1,1+4);
	\draw[thick] (2,1+4) -- (2,0+4) -- (3,0+4) -- (3,1+4);
	\draw[thick] (2+2,1+4) -- (2+2,0+4) -- (3+2,0+4) -- (3+2,1+4);
	\draw[thick] (2+4,1+4) -- (2+4,0+4) -- (3+4,0+4) -- (3+4,1+4);
	\draw[thick] (2+6,1+4) -- (2+6,0+4) -- (3+6,0+4) -- (3+6,1+4);
	\draw[thick] (2+8,1+4) -- (2+8,0+4) -- (3+8,0+4) -- (3+8,1+4);
	\draw[thick] (2+10,1+4) -- (2+10,0+4) -- (3+10,0+4) -- (3+10,1+4);
	
	\node[circle,shading=ball,minimum width=0.4cm,color=black!20!white] (ball) at (0.3,4.2) {};
	\node[circle,shading=ball,minimum width=0.4cm,color=black!20!white] (ball) at (0.7+2,4.2) {};
	\node[circle,shading=ball,minimum width=0.4cm,color=black!20!white] (ball) at (0.3+10,4.2) {};
	\node[circle,shading=ball,minimum width=0.4cm,color=black!20!white] (ball) at (0.7+12,4.2) {};
	
	\node[circle,shading=ball,minimum width=0.4cm,color=black!20!white] (ball) at (0.3,2.2) {};
	\node[circle,shading=ball,minimum width=0.4cm,color=black!20!white] (ball) at (0.7+2,2.2) {};
	\node[circle,shading=ball,minimum width=0.4cm,color=black!20!white] (ball) at (0.3+10,2.2) {};
	\node[circle,shading=ball,minimum width=0.4cm,color=black!20!white] (ball) at (0.75,2.2) {};
	\node[circle,shading=ball,minimum width=0.4cm,color=black!20!white] (ball) at (0.7+6,2.2) {};
	\node[circle,shading=ball,minimum width=0.4cm,color=black!20!white] (ball) at (0.7+12,2.2) {};
	\node[circle,shading=ball,minimum width=0.4cm,color=black!20!white] (ball) at (0.25+12,2.2) {};
	\node[circle,shading=ball,minimum width=0.4cm,color=black!20!white] (ball) at (0.75+10,2.2) {};
	
	\node[circle,shading=ball,minimum width=0.4cm,color=black!20!white] (ball) at (0.3,0.2) {};
	\node[circle,shading=ball,minimum width=0.4cm,color=black!20!white] (ball) at (0.7+2,0.2) {};
	\node[circle,shading=ball,minimum width=0.4cm,color=black!20!white] (ball) at (0.3+10,0.2) {};
	\node[circle,shading=ball,minimum width=0.4cm,color=black!20!white] (ball) at (0.75,0.2) {};
	\node[circle,shading=ball,minimum width=0.4cm,color=black!20!white] (ball) at (0.7+6,0.2) {};
	\node[circle,shading=ball,minimum width=0.4cm,color=black!20!white] (ball) at (0.7+12,0.2) {};
	\node[circle,shading=ball,minimum width=0.4cm,color=black!20!white] (ball) at (0.75+10,0.2) {};
	\node[circle,shading=ball,minimum width=0.4cm,color=black!20!white] (ball) at (0.525,0.55) {};
	\node[circle,shading=ball,minimum width=0.4cm,color=black!20!white] (ball) at (0.25+2,0.2) {};
	\node[circle,shading=ball,minimum width=0.4cm,color=black!20!white] (ball) at (0.25+12,0.2) {};
	\node[circle,shading=ball,minimum width=0.4cm,color=black!20!white] (ball) at (0.475+12,0.55) {};
	\node[circle,shading=ball,minimum width=0.4cm,color=black!20!white] (ball) at (0.25+6,0.2) {};
	
	\draw[->,thick] (-0.2,4.3) to[bend right] (-0.2,2.7);
	\draw[->,thick] (-0.2,2.3) to[bend right] (-0.2,0.7);
\end{tikzpicture}
\caption{Illustration of the cell interpretation with $n=7$ and $s=4$ in the first three drawings.}
\label{fig1}
 \end{figure}

The coupon collector's problem is an old problem of probability theory which in its simplest form dates back to de Moivre, Laplace and Euler, see  \cite{Euler},  \cite{Laplace} and  \cite{DeMoivre}.  Whereas de Moivre used a die with $ s $ faces to pose the problem, Euler and Laplace used a lottery interpretation as motivation. However, a more recent example for a situation in which the  coupon collector's problem occurs is the collection of pictures of the participating players of all teams before and during every World Cup. Typically, fans can buy the pictures in packages of five or six. Two natural questions which arise are:  How many packages does one need to buy to get the full or a specific portion of the full set of players? How many stickers are missing after one has bought $ k $ packages? The first  question was for example studied in \cite{BaumBillingsley}, \cite{Johnson}, \cite{Kolchin} or \cite{Stadje}. In the work at hand we will deal with the latter of the two problems. 
 The version  of the coupon collector's problem we consider can be described as follows. Assume there are $ n $ different coupons. At each time we draw $ s $ of these $ n $ coupons, where we assume that each of the $ {n \choose s} $ choices occurs with the same probability.
 We are then interested in the distribution of the number $ Z_{n,s}(k_n) $ of coupons that have not been drawn in the first $ k=k_n $ drawings. In a conceptually equivalent interpretation the $ n $ coupons are represented by $ n $ different cells numbered $ 1, \ldots, n $, and in each drawing we place $ s $ particles into $ s $ distinct cells, see Figure \ref{fig1}. We are then interested in the distribution of the number $ Z_{n,s}(k_n) $ of  empty cells after $ k=k_n$ drawings.
 
 The behaviour of $ Z_{n,s}(k_n) $ for the different regimes of $ k_n $ has been subject to numerous research works over the years.
In \cite{Mahmoud}  convergence towards a normal limit is proved for the sublinear and linear regime, that is, for $k_n= o(n)$ and $k_n=\alpha n $, respectively. In \cite{Smythe} the author proves a central limit theorem for a generalized coupon collector's problem allowing for \textit{random} package sizes $ S $ in the lower superlinear regime, i.e., for $\frac{k_n}{n}\rightarrow \infty \text{ and } \limsup_n \frac{k_n}{n \log(n)}< \frac{1}{\EE[S]}$ using a martingale representation.
In  \cite{MikhailovNormal80} the method of moments is used to prove normal approximation for the case, where $n(\frac{k_n s}{n})^r \rightarrow \infty  $ for all $ r  \in  \NN $ and $\EE[Z_{n,s}(k_n) ] \rightarrow \infty$, which covers our regime of normal approximation introduced below. However, no rates of convergence are given in neither of the works mentioned so far. In \cite{VatutinMikhailov} the authors deduce rates of convergence in the Kolmogorov distance towards a normal limit for $ k_n=o(n \log(n))$ which are of order $1/\sqrt{\var(Z_{n,s}(k_n))}$. We complement these bounds in the case where $ k_n $ is assumed to be of the form $ k_n=\frac{n}{s}(\alpha \log(n)+x) $ for some $ \alpha \in (0,1 )$ and $ x \in \RR $. In the case $ \alpha=1 $, i.e.\ if $k_n =\frac{n}{s}( \log(n)+x)$, the author in \cite{MikhailovPoisson77} uses the Stein-Chen method to prove convergence towards a Poisson limit for the number of cells containing exactly $ r $ particles in a more general setting allowing for multiple particles being placed in one cell at each step. The same question is also studied in \cite{SchillingHenze} and a Poisson limit is deduced. In the same work, the authors show that for the special case $ r=0 $ the condition of equally probable group drawings can be relaxed to a certain extend without losing the limiting Poisson distribution. In \cite[Theorem 6.F]{BHS92} the authors also prove rates of convergence of order $\frac{\log(n)}{n}$ towards a Poisson distribution in this regime, which are optimal in view of \cite[Theorem 3.D]{BHS92} . For $ r=0 $ this covers our setting of Poisson approximation with the same rates, which are included here only for completeness and to demonstrate that both limit theorems can be based on the same coupling argument.

 As we have explained above, there exists a sharp asymptotic distributional phase transition at $ \alpha=1 $ in the sense that for $ \alpha \in (0,1) $ the random variable  $ Z_{n,s}(k_n) $ asymptotically follows a normal distribution whereas for $ \alpha=1 $ we obtain a Poisson limit. However, in both cases we use Stein's method in combination with the same size-biased coupling construction to prove upper bounds on the distance between $ Z_{n,s}(k_n) $ and Gaussian or Poisson random variable, respectively. Our results are presented in the next section, while the coupling construction is explained in Section \ref{sec:Coupling}. The proof of the normal approximation result for $\alpha\in(0,1)$ is the content of Section \ref{sec:NormalProof}, while the Poisson limit theorem is derived in Section \ref{sec:PoissonProof}.

\section{Results}

Denote by $ (C_i)_{i=1,\ldots n} $ the collection of cells in the coupon collector's problem and let $ Z_{n,s}(k_n) $ be the number of empty cells after $ k_n $ drawings, i.e.
\[
Z_{n,s}(k_n) = \sum_{j=1}^{n} \mathbf{1}_{E_{n,j}(k_n)},
\]
where  $ E_{n,j}(k)= \{|C_j|=0\text{ after } k_n \text{ drawings}\} $ and where $|C_j|$ stands for the number of particles in cell $C_j$. Moreover, we shall denote for two random variables $X$ and $Y$ by
$$
d_W(X,Y) := \sup_{h\in{\rm Lip}(1)}\big|\EE[h(X)]-\EE[h(Y)]\big|
$$
the Wasserstein distance between $X$ and $Y$, where the supremum runs over all Lipschitz functions $h:\RR\to\RR$ with Lipschitz constant less than or equal to one.  
We consider the case where  $ k_n=\frac{n}{s}(\alpha\log(n)+x) $ for fixed $ s \in \NN $, $ x \in \RR $ and $ \alpha \in (0,1) $. Note that since $ k_n $ denotes the number of drawings we always assume implicitly that $ k_n$ is an integer, in particular we assume that $ n $ is large enough so that $ k_n \geq 0 $. Furthermore, we define the centred and normalized random variables
\begin{align}\label{eq:StandardizedZ}
\tilde{Z}_{n,s}(k_n):= \frac{Z_{n,s}(k_n)- \EE[Z_{n,s}(k_n)]}{\sqrt{\var(Z_{n,s}(k_n))}}.
\end{align}
Throughout the paper we use the notation $ C(x_1, x_2,\ldots)  $ to indicate that a constant $C\in(0,\infty)$ only depends on parameters $ x_1, x_2,\ldots $ of the model.

\begin{theorem}\label{Thm:Normal}
	Put $k_n= \frac{n}{s}(\alpha\log(n)+x) $ for some $s\in\NN$, $ x \in \RR $ and $ \alpha \in (0,1) $. Let $ Z_{n,s}(k_n) $ be the number of empty cells after $ k_n $ drawings as introduced above and denote by $ G $ a standard Gaussian random variable. Then there exist constants $ C=C(s,x) \in (0, \infty) $ and $ N=N(s,x,\alpha) \in \NN $ such that for all $ n \geq N $,
	\begin{align}\label{eq:RateNormal}
	d_{W}(\tilde{Z}_{n,s}(k_n), G) \leq \begin{cases}
		C\frac{\sqrt{\log(n)}}{n^{\alpha}} &\text{ for } \alpha \in (0, \frac{1}{3}]\\[0.2cm]
		Cn^{-{1\over 2}(1-\alpha)}&\text{ for } \alpha \in( \frac{1}{3},1).
	\end{cases}
	\end{align}
\end{theorem}

\begin{remark}\rm 
\begin{itemize}
	\item[(i)] Since $n^{-\frac{1}{2}(1-\alpha)}$ is of the same order as $1/\sqrt{\var(Z_{n,s}(k_n))}$ we believe that the rate in Theorem \ref{Thm:Normal} is optimal at least in the regime $\alpha\in(1/3,1)$. We leave it as an open problem to decide whether or not the rate is optimal for $\alpha\in(0,1/3]$. On the other hand, Remark \ref{rm:OptimalVar} shows that in this situation the rate cannot be improved on the basis of the general normal approximation bound \eqref{NormalBound} we use.
	\item[(ii)] One might ask whether the Wasserstein distance can be replaced by the Kolmogorov distance $d_{K}(\tilde{Z}_{n,s}(k_n), N)=\sup_{u\in\RR}|\PP(\tilde{Z}_{n,s}(k_n)\leq u)-\PP(N\leq u)|$. As we shall explain in Remark \ref{rem:Kolmogorov} in more detail, this is not possible by means of the size-biased coupling approach of Stein's method for normal approximation using our coupling construction. In fact, the resulting bound does in this case not even tend to zero with $n$. On the other hand, the presumably suboptimal bound $d_{K}(\tilde{Z}_{n,s}(k_n), N)\leq C(\log(n))^{1/4}n^{-\alpha/2}$ for $0<\alpha\leq 1/3$ and $d_{K}(\tilde{Z}_{n,s}(k_n), N)\leq Cn^{-\frac{1}{4}(1-\alpha)}$ for $1/3<\alpha<1$ follows directly from the fact that the Kolmogorov distance can always be bounded by the square-root of the Wasserstein distance. 
\end{itemize}
\end{remark}

The next result complements Theorem \ref{Thm:Normal} by considering the case $\alpha=1$ for which the upper bound \eqref{eq:RateNormal} does not tend to zero, as $n\to\infty$. As emphasized already above, the result is known from \cite{BHS92} and is included here only for completeness. Similarly as above, for two random variables $X$ and $Y$ we denote by
$$
d_{TV}(X,Y) := \sup_A\big|\PP(X\in A)-\PP(Y\in B)\big|
$$
the total variation distance between $X$ and $Y$, where the supremum is taken over all Borel sets $A\subset\RR$.

\begin{theorem}\label{Thm:Poisson}
	Put $ k_n= \frac{n}{s}(\log(n)+x) $ for some $s\in\NN$ and $ x \in \RR $. Let $ Z_{n,s}(k_n) $ be the number of empty cells after $ k_n $ drawings as introduced above and denote by $ W $ a Poisson random variable with parameter $ \lambda_n=\EE[Z_{n,s}(k_n) ]$. Then there exist constants $ \tilde{C}=\tilde{C}(s,x) \in (0, \infty) $ and $\tilde{N}=\tilde{N}(s,x) \in \NN $ such that for all $ n \geq \tilde{N} $,
	\begin{align*}
	d_{TV}(Z_{n,s}(k_n), W) \leq  \tilde{C}\, \frac{ \log(n)}{n}.
	\end{align*}
\end{theorem}

For the proof of both results we use Stein's method in combination with a size-biased coupling. We start by describing this coupling in the next section.

\section{Coupling construction}\label{sec:Coupling}

 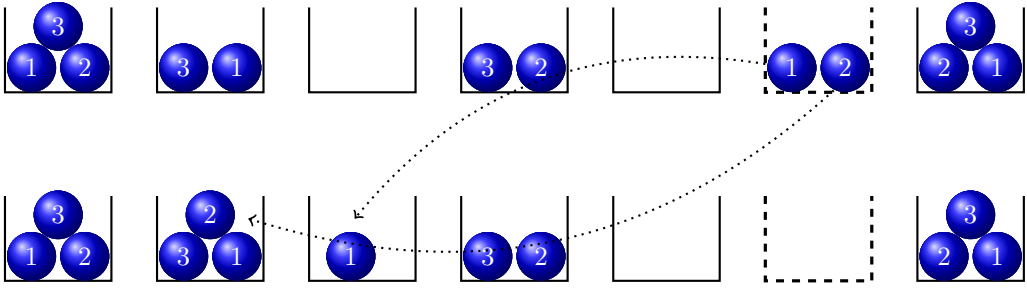
\begin{figure}[t]
	\centering
	\begin{tikzpicture}
		\draw[thick] (-0.2,1) -- (-0.2,-0.125) -- (1.2,-0.125) -- (1.2,1);
		\draw[thick] (2-0.2,1) -- (2-0.2,-0.125) -- (3+0.2,-0.125) -- (3+0.2,1);
		\draw[thick] (2+2-0.2,1) -- (2+2-0.2,-0.125) -- (3+2+0.2,-0.125) -- (3+2+0.2,1);
		\draw[thick] (2+4-0.2,1) -- (2+4-0.2,-0.125) -- (3+4+0.2,-0.125) -- (3+4+0.2,1);
		\draw[thick] (2+6-0.2,1) -- (2+6-0.2,-0.125) -- (3+6+0.2,-0.125) -- (3+6+0.2,1);
		\draw[very thick,dashed] (2+8-0.2,1) -- (2+8-0.2,-0.125) -- (3+8+0.2,-0.125) -- (3+8+0.2,1);
		\draw[thick] (2+10-0.2,1) -- (2+10-0.2,-0.125) -- (3+10+0.2,-0.125) -- (3+10+0.2,1);

		\node[circle,shading=ball,minimum width=0.4cm,color=black!20!white] (ball) at (0.15,0.2) {\small\textcolor{white}{$1$}};
		\node[circle,shading=ball,minimum width=0.4cm,color=black!20!white] (ball) at (0.85+2,0.2) {\small\textcolor{white}{$1$}};
		\node[circle,shading=ball,minimum width=0.4cm,color=black!20!white] (ball) at (0.15+10,0.2) {\small\textcolor{white}{$1$}};
		\node[circle,shading=ball,minimum width=0.4cm,color=black!20!white] (ball) at (0.85,0.2) {\small\textcolor{white}{$2$}};
		\node[circle,shading=ball,minimum width=0.4cm,color=black!20!white] (ball) at (0.85+6,0.2) {\small\textcolor{white}{$2$}};
		\node[circle,shading=ball,minimum width=0.4cm,color=black!20!white] (ball) at (0.85+12,0.2) {\small\textcolor{white}{$1$}};
		\node[circle,shading=ball,minimum width=0.4cm,color=black!20!white] (ball) at (0.85+10,0.2) {\small\textcolor{white}{$2$}};
		\node[circle,shading=ball,minimum width=0.4cm,color=black!20!white] (ball) at (0.5,0.75) {\small\textcolor{white}{$3$}};
		\node[circle,shading=ball,minimum width=0.4cm,color=black!20!white] (ball) at (0.15+2,0.2) {\small\textcolor{white}{$3$}};
		\node[circle,shading=ball,minimum width=0.4cm,color=black!20!white] (ball) at (0.15+12,0.2) {\small\textcolor{white}{$2$}};
		\node[circle,shading=ball,minimum width=0.4cm,color=black!20!white] (ball) at (0.5+12,0.75) {\small\textcolor{white}{$3$}};
		\node[circle,shading=ball,minimum width=0.4cm,color=black!20!white] (ball) at (0.15+6,0.2) {\small\textcolor{white}{$3$}};
		
	   	\draw[thick] (-0.2,1-2.5) -- (-0.2,-0.125-2.5) -- (1.2,-0.125-2.5) -- (1.2,1-2.5);
	   \draw[thick] (2-0.2,1-2.5) -- (2-0.2,-0.125-2.5) -- (3+0.2,-0.125-2.5) -- (3+0.2,1-2.5);
	   \draw[thick] (2+2-0.2,1-2.5) -- (2+2-0.2,-0.125-2.5) -- (3+2+0.2,-0.125-2.5) -- (3+2+0.2,1-2.5);
	   \draw[thick] (2+4-0.2,1-2.5) -- (2+4-0.2,-0.125-2.5) -- (3+4+0.2,-0.125-2.5) -- (3+4+0.2,1-2.5);
	   \draw[thick] (2+6-0.2,1-2.5) -- (2+6-0.2,-0.125-2.5) -- (3+6+0.2,-0.125-2.5) -- (3+6+0.2,1-2.5);
	   \draw[very thick, dashed] (2+8-0.2,1-2.5) -- (2+8-0.2,-0.125-2.5) -- (3+8+0.2,-0.125-2.5) -- (3+8+0.2,1-2.5);
	   \draw[ thick] (2+10-0.2,1-2.5) -- (2+10-0.2,-0.125-2.5) -- (3+10+0.2,-0.125-2.5) -- (3+10+0.2,1-2.5);
	   
	   \node[circle,shading=ball,minimum width=0.4cm,color=black!20!white] (ball) at (0.15,0.2-2.5) {\small\textcolor{white}{$1$}};
	   \node[circle,shading=ball,minimum width=0.4cm,color=black!20!white] (ball) at (0.85+2,0.2-2.5) {\small\textcolor{white}{$1$}};
	   \node[circle,shading=ball,minimum width=0.4cm,color=black!20!white] (ball) at (0.35+4,0.2-2.5) {\small\textcolor{white}{$1$}};
	   \node[circle,shading=ball,minimum width=0.4cm,color=black!20!white] (ball) at (0.85,0.2-2.5) {\small\textcolor{white}{$2$}};
	   \node[circle,shading=ball,minimum width=0.4cm,color=black!20!white] (ball) at (0.85+6,0.2-2.5) {\small\textcolor{white}{$2$}};
	   \node[circle,shading=ball,minimum width=0.4cm,color=black!20!white] (ball) at (0.85+12,0.2-2.5) {\small\textcolor{white}{$1$}};
	   \node[circle,shading=ball,minimum width=0.4cm,color=black!20!white] (ball) at (0.5+2,0.75-2.5) {\small\textcolor{white}{$2$}};
	   \node[circle,shading=ball,minimum width=0.4cm,color=black!20!white] (ball) at (0.5,0.75-2.5) {\small\textcolor{white}{$3$}};
	   \node[circle,shading=ball,minimum width=0.4cm,color=black!20!white] (ball) at (0.15+2,0.2-2.5) {\small\textcolor{white}{$3$}};
	   \node[circle,shading=ball,minimum width=0.4cm,color=black!20!white] (ball) at (0.15+12,0.2-2.5) {\small\textcolor{white}{$2$}};
	   \node[circle,shading=ball,minimum width=0.4cm,color=black!20!white] (ball) at (0.5+12,0.75-2.5) {\small\textcolor{white}{$3$}};
	   \node[circle,shading=ball,minimum width=0.4cm,color=black!20!white] (ball) at (0.15+6,0.2-2.5) {\small\textcolor{white}{$3$}};
		
		\draw[->,thick,dotted] (9.8,0.25) to[bend right] (4.4,-1.8);
		\draw[->,thick,dotted] (10.7,-0.1) to[bend left] (3,-1.8);
	\end{tikzpicture}
	\caption{Illustration of the coupling construction (continuation of Figure \ref{fig1}); the artificially isolated cell with label $I$ is the dashed one.}
	\label{fig2}
\end{figure}

\noindent We are dealing with the coupon collector's problem with $ n $ coupons, which can be interpreted as $ n $ distinct cells. At each time step we place a fixed number $ s \leq n $ of particles in $ s $ different cells. In order to keep track of when particles are placed into cells, we label all particles in the $ m $th drawing with the letter $ m $. Now, after $ k_n $ placements, we choose one of the $ n $ cells, which we denote by $ C_I $,  uniformly at random and take all particles out of it. For a particle labelled $ j  $ taken from $ C_I $ we now choose one of the $n-s $ cells not containing a particle with label $ j  $ uniformly and place the particle into it. We proceed in the same manner until all particles from cell $ I $ have been redistributed into the remaining $ n-1 $ cells, see Figure \ref{fig2}. Denote by $ F_{n,j}^I $ the event that at least one particle from cell $ I $ is placed into cell $ j $ and  put $\bar{F}_{n,j}^I:=(F_{n,j}^I)^c$. Furthermore, we define $ E^I_{n,j}(k_n): =  E_{n,j}(k_n) \cap \bar{F}_{n,j}^I.   $ For any $ j \neq I $ we then have
\begin{align}\label{eq:IndicatorSizeBias}
\PP(E^I_{n,j}(k_n))&=\sum_{\ell=0}^{k_n} \PP(E_{n,j}(k_n)\cap \bar{F}_{n,j}^I\cap\{|C_I|=\ell\})\notag\\
&= \sum_{\ell=0}^{k_n} \PP( \bar{F}_{n,j}^I\,\vert\, E_{n,j}(k_n)\cap\{ |C_I|=\ell\}) \  \PP( |C_I|=\ell\,|\,E_{n,j}(k_n))\ \PP(E_{n,j}(k_n)).
\end{align}
Note that by construction we have  
\begin{align}
\PP( E_{n,j}(k_n))&= \Big(1-\frac{s}{n}\Big)^{k_n}, \label{ProbEj}\\
\PP( |C_I|=\ell\,|\,E_{n,j}(k_n)) &= {k_n \choose \ell } \Big(\frac{s}{n-1}\Big)^\ell \Big(1-\frac{s}{n-1}\Big)^{k_n-\ell}, \label{PropCIEj}
\end{align}
since conditioning on the event $ E_{n,j}(k_n) $ simply means that we can only place particles in $ n-1 $ instead of $n$ cells. In addition,
\begin{align*}
\PP( \bar{F}_{n,j}^I\,\vert\, E_{n,j}(k_n)\cap\{ |C_I|=\ell\}) &=\Big(1-\frac{1}{n-s}\Big)^\ell,
\end{align*}
since for each of the particles from cell $I$ the probability that it is placed into an originally empty cell is $ \frac{1}{n-s} $.
Putting these expressions back into \eqref{eq:IndicatorSizeBias} we obtain
\begin{align}\label{SBIndicatorProb}
\PP(E_{n,j}^{I}(k_n))&= \sum_{\ell=0}^{k_n}  {k_n \choose \ell } \Big(\frac{s}{n-1}\Big)^\ell \Big(1-\frac{s}{n-1}\Big)^{k_n-\ell}\Big(1-\frac{1}{n-s}\Big)^\ell \, \Big(1-\frac{s}{n}\Big)^{k_n} \notag\\[0.2cm]
&=\Big(1-\frac{s}{n-1}\Big)^{k_n}\notag\\[0.2cm]
&
 =\PP(E_{n,j}(k_n)\,\vert\, E_{n,i}(k_n))
 \end{align}
 for any $ i \neq j $. We thus conclude that the random variable
 \begin{align*}
 Z^I_{n,s}(k_n) := 1+\sum_{\substack{j=1\\ j \neq I}}^{n} \mathbf{1}_{E^I_{n,j}(k_n)}
 \end{align*}
has the $ Z_{n,s}(k_n)- $size biased distribution, meaning that
$$
\PP( Z^I_{n,s}(k_n) = y) = {y\over \EE[Z_{n,s}(k_n)]}\PP(Z_{n,s}(k_n) = y),\qquad y\in\{0,1,\ldots,n\},
$$
see \cite{ArratiaGoldsteinKochman}.

\begin{remark}\label{rem:Kolmogorov}\rm
	Theorem 5.6 in \cite{ArratiaGoldsteinKochman} provides a bound on the Kolmogorov distance between $\tilde{Z}_{n,s}(k_n)$ and a standard Gaussian random variable using a size-biased coupling. However, for this to yield a central limit theorem one needs that   $ |Z^I_{n,s}-Z_{n,s}|  =o(\sqrt{n})$ almost surely, as $n\to\infty$. For the coupling described above we have $ |Z^I_{n,s}-Z_{n,s}| =n-s-1 $ if in all $ k_n $ drawings the same $ s $ cells are filled and the remaining $ n-s $ cells are filled when redistributing the $ k_n $ particles of one of the filled cells. Consequently, the result in \cite{ArratiaGoldsteinKochman} does not lead to a meaningful bound on the Komogorov distance as already explained in part (ii) of the remark after Theorem \ref{Thm:Normal}.
\end{remark}

\section{Proof of Theorem~\ref{Thm:Normal}}\label{sec:NormalProof}

Following \cite[Theorem 1.1]{GR96} the Wasserstein distance between $\tilde{Z}_{n,s}=\tilde{Z}_{n,s}(k_n)$ as defined in \eqref{eq:StandardizedZ} and  a standard Gaussian random variable $G$ is bounded by
\begin{align}\label{NormalBound}
d_W(\tilde{Z}_{n,s},G)\leq \frac{\lambda_n}{\sigma_n^2}\sqrt{\var(\EE[Z_{n,s}-Z^I_{n,s}\vert \cC_n(k_n)])}+ \frac{\lambda_n}{\sigma_n^3} \EE[(Z_{n,s}-Z^I_{n,s})^2],
\end{align}
where $ \cC_n(k_n) $ denotes the configuration of the $ n $ cells after $ k_n $ drawings, $ \lambda_n=\EE[Z_{n,s}(k_n)] $ and $ \sigma_n^2= \var(Z_{n,s}(k_n)) $. In the next sections we bound further the right-hand side of \eqref{NormalBound} by dealing with the individual terms.

\subsection{Expectation and variance}

We start by bounding from above the expectation and from below the variance of $ Z_{n,s}$.  First, we note that  for $ k_n=\frac{n}{s}(\alpha\log(n)+x) $ 
 we have
\begin{align}\label{eq:ProbEmptyCellAlpha}
\PP( E_{n,j}(k_n))&=\Big(1-\frac{s}{n}\Big)^{k_n} =\bigg(\Big(1-\frac{s}{n}\Big)^{n/s}\bigg)^{\alpha\log(n)+x} \sim \ e^{-(\alpha\log(n)+x) }=\frac{\, e^{-x}}{n^\alpha},
\end{align}
where we write $ f(n) \sim g(n)  $ for two functions $f,g:\NN\to\RR$ if $\displaystyle{\lim_{n\rightarrow \infty} \frac{f(n)}{g(n)}=1}$. 
Thus,
\begin{align}\label{Expectation}
\EE[Z_{n,s}(k_n)]= \sum_{j=1}^{n}\PP( E_{n,j}(k_n)) \sim e^{-x}n^{1-\alpha}.
\end{align}
In particular, we can find $ n_0 \in \NN $, depending on $ \alpha $, such that for all $ n \geq n_0 $,
\begin{align} \label{eq:UpperBoundPEnj}
\PP( E_{n,j}(k_n))&\leq \frac{2\, e^{-x}}{n^\alpha}
\end{align}

Similarly, we obtain
\begin{align}\label{eq:ProbEmptyCellsDiffAlpha}
\Big(1-\frac{s}{n}\Big)^{k_n} -\Big(1-\frac{s}{n-1}\Big)^{k_n}&=\Big(1-\frac{s}{n}\Big)^{k_n}\bigg(1-\Big(\frac{(n-1-s)n}{(n-1)(n-s)}\Big)^{k_n} \bigg)\notag \\
&\leq\Big(1-\frac{s}{n}\Big)^{k_n}\bigg(-k_n\Big(\frac{(n-1-s)n}{(n-1)(n-s)}-1\Big)\bigg)\notag\\
&=\Big(1-\frac{s}{n}\Big)^{k_n} \, k_n\, \frac{(n-1)(n-s)-(n-1-s)n}{(n-1)(n-s)}\notag\\
&=\Big(1-\frac{s}{n}\Big)^{k_n} \frac{k_n\,s}{(n-1)(n-s)}\notag\\
&\leq \frac{2 e^{-x} (\alpha \log(n)+x)}{n^{1+\alpha}},
\end{align}
for all $ n \geq n_1 $ for some $ n_1= n_1(\alpha,s) \in \NN $,
where we have used that 
\begin{align}\label{Ineq:1-z^k}
1-z^k\leq-k(z-1) 
\end{align}
for   $ z \in (0,1) $.
Using that $\var({\bf 1}_A)=\PP(A)(1-\PP(A))$ and $\Cov({\bf 1}_A,{\bf 1}_B)=\PP(A\cap B)-\PP(A)\PP(B)=\PP(B)(\PP(A|B)-\PP(A))$, we see that
\begin{align*}
	\var(Z_{n,s})&= \sum_{j=1}^{n} \var[\mathbf{1}_{E_{n,j}(k_n)}]+ \sum_{j=1}^{n} \sum_{\substack{i=1 \\ i\neq j}}^n \Cov[\mathbf{1}_{E_{n,j}(k_n)}, \mathbf{1}_{E_{n,i}(k_n)}]\\
	&=\sum_{j=1}^n\PP(E_{n,j}(k_n))(1-\PP(E_{n,j}(k_n)))\\
	&\qquad\qquad\qquad+\sum_{j=1}^{n} \sum_{\substack{i=1 \\ i\neq j}}^n[\PP(E_{n,i}(k_n))(\PP(E_{n,j}(k_n)\,|\,E_{n,i}(k_n))-\PP(E_{n,j}(k_n)))]\\
	&=\bigg(1-\Big(1-\frac{s}{n}\Big)^{k_n}+(n-1)\Big( \Big(1-\frac{s}{n-1}\Big)^{k_n} -\Big(1-\frac{s}{n}\Big)^{k_n} \Big)\bigg )\EE[Z_{n,s}].
\end{align*}
Applying now \eqref{eq:ProbEmptyCellsDiffAlpha} we conclude that there exists a constant $ c_1(s,x) \in (0, \infty) $ such that for $  n \geq \max\{n_0,n_1\} $ the lower variance bound
\begin{align}\label{Variance}
\var(Z_{n,s})&\geq\Big(1-\frac{2 e^{-x}}{n^\alpha}-\frac{2e^{-x}\, (\alpha \log(n)+x)}{ n^\alpha}\Big)\EE[Z_{n,s}]  \geq c_1(s,x)\,  \EE[Z_{n,s}] ,
\end{align}
holds.


\subsection{Bounding $\var\big(\EE[Z_{n,s}-Z^I_{n,s}\vert\cC_n(k_n) ]\big)$}
For the variance of the conditional expectation on the right-hand side in \eqref{NormalBound} we obtain
\begin{align}\label{VarianceTerms}
&\var\big(\EE[Z_{n,s}-Z^I_{n,s}\vert\cC_n(k_n) ]\big)=\var\Big(\EE\Big[\sum_{j=1}^n \mathbf{1}_{E_{n,j}(k_n)\cap F_{n,j}^I}\vert\cC_n(k_n) \Big]\Big)
\notag\\
&=\EE\Big[\Big(\sum_{j=1}^n\EE[\mathbf{1}_{E_{n,j}(k_n)\cap F_{n,j}^I}\vert \cC_n(k_n) ]\Big)^2\Big]- \EE\Big[\sum_{j=1}^n \EE[\mathbf{1}_{E_{n,j}(k_n)\cap F_{n,j}^I}\vert \cC_n(k_n) ]\Big]^2 \notag\\
&=\EE\Big[\sum_{j=1}^n \EE[\mathbf{1}_{E_{n,j}(k_n)\cap F_{n,j}^I}\vert \cC_n(k_n) ]^2\Big]\notag\\
&\quad+\Bigg(\EE\Big[\sum_{j=1}^n \sum_{\substack{i=1 \\ i\neq j}}^n \,\EE[\mathbf{1}_{E_{n,j}(k_n)\cap F_{n,j}^I}\vert \cC_n(k_n) ]\, \EE[\mathbf{1}_{E_{n,i}(k_n)\cap F_{n,i}^I}\vert \cC_n(k_n) ]\Big]  - \Big[\sum_{j=1}^n \EE[\mathbf{1}_{E_{n,j}(k_n)\cap F_{n,j}^I}]\Big]^2  \Bigg)\notag\\
&=: T_1+T_2.
\end{align}
We start by dealing with $ T_1 $. Note that conditionally on the event $ E_{n,j}(k_n) $ we have
\begin{equation}\label{ConditionalExpectation}
\EE[\mathbf{1}_{F_{n,j}^I}\vert \cC_n(k_n) ]= 1-\Big(1-\frac{1}{n-s}\Big)^{|C_I|},
\end{equation}
so that with \eqref{Ineq:1-z^k} we obtain
\begin{align*}
T_1= \EE\Big[\sum_{j=1}^n \mathbf{1}_{E_{n,j}(k_n)} \Big(1-\Big(1-\frac{1}{n-s}\Big)^{|C_I|}\Big)^2\Big]&\leq \frac{1}{(n-s)^2}\EE\Big[\sum_{j=1}^n \mathbf{1}_{E_{n,j}(k_n)}|C_I|^2\Big].
\end{align*}
 Denoting by $ D^I_m $ the event that a particle is placed into cell $ C_I $ in the $m$th  drawing for some $m\in\{1,\ldots,k_n\}$ we see that
\begin{align*}
\EE[\mathbf{1}_{E_{n,j}(k_n)}|C_I|^2]&=  \EE\Big[\mathbf{1}_{E_{n,j}(k_n)}\Big(\sum_{m=1}^{k_n}\mathbf{1}_{D_m^I} \Big)^2 \Big]\\
& =  \EE\Big[\sum_{m=1}^{k_n}\mathbf{1}_{E_{n,j}(k_n)}\mathbf{1}_{D_m^I}  \Big] + \EE\Big[\sum_{m=1}^{k_n} \sum_{\substack{r=1 \\ r\neq m}}^{k_n}\mathbf{1}_{E_{n,j}(k_n)}\mathbf{1}_{D_m^I}\mathbf{1}_{D_r^I}  \Big]\\
& = \PP(E_{n,j}(k_n)) \Big(\sum_{m=1}^{k_n}\PP(D_m^I  \big\vert E_{n,j}(k_n)) + \sum_{m=1}^{k_n} \sum_{\substack{r=1 \\ r\neq m}}^{k_n}\PP(D_m^I \cap D_r^I \big\vert E_{n,j}(k_n)) \Big).
\end{align*}
Since all $s$-placements occur with the same probability and conditioning on the event $ E_{n,j}(k_n) $ simply means that we can only place particles into $ n-1 $ cells,  we have
\[
\PP(D_m^I  \big\vert E_{n,j}(k_n))= \frac{s}{n-1}.
\]
As the drawings are independent of each other we get
\begin{align*}
&\sum_{m=1}^{k_n}\PP(D_m^I  \big\vert E_{n,j}(k_n)) + \sum_{m=1}^{k_n} \sum_{\substack{r=1 \\ r\neq m}}^{k_n}\PP(D_m^I \cap D_r^I \big\vert E_{n,j}(k_n))\\
&=  \sum_{m=1}^{k_n}\PP(D_m^I  \big\vert E_{n,j}(k_n)) + \sum_{m=1}^{k_n} \sum_{\substack{r=1 \\ r\neq m}}^{k_n}\PP(D_m^I  \big\vert E_{n,j}(k_n)) \PP(D_r^I \big\vert E_{n,j}(k_n))=\frac{k_n s}{n-1} + \frac{k_n (k_{n}-1)s^2}{(n-1)^2}.
\end{align*}
Using also \eqref{eq:ProbEmptyCellAlpha} we obtain that there exists 
a constant  $c_2(s,x) \in (0,\infty)$ such that we can bound  $T_1$ by
\begin{align}\label{T1}
T_1 \leq \Big(1-\frac{s}{n}\Big)^{k_n}\frac{n}{(n-s)^2}\Big(\frac{k_n s}{n-1} + \frac{k_n (k_{n}-1)s^2}{(n-1)^2}\Big)\leq  c_2(s,x) \frac{\alpha^2 \log(n)^2}{n^{1+\alpha}}.
\end{align}
To deal with the term $ T_2 $ first note that by \eqref{ProbEj} and \eqref{PropCIEj},
\begin{align}\label{ProbEjFj}
 \EE[\mathbf{1}_{E_{n,j}(k_n)\cap F_{n,j}^I}]&= \PP(E_{n,j}(k_n)) \PP( F_{n,j}^I \vert E_{n,j}(k_n)) \notag \\
 &= \PP(E_{n,j}(k_n)) \sum_{\ell=0}^{k_n} \PP( F_{n,j}^I \vert E_{n,j}(k_n) \cap \{\vert C_I \vert = \ell\} ) \PP( C_I \vert = \ell \vert  E_{n,j}(k_n))\notag \\
 &=\Big(1-\frac{s}{n}\Big)^{k_n} \sum_{\ell=0}^{k_n} \Big(1-\big(1-\frac{1}{n-s}\big)^\ell\Big) {k_n \choose \ell} \Big(\frac{s}{n-1}\Big)^\ell \Big(1-\frac{s}{n-1}\Big)^{k_n -\ell}\notag \\
&= \Big(1-\frac{s}{n}\Big)^{k_n} \Big(1- \sum_{\ell=0}^{k_n}\big(1-\frac{1}{n-s}\big)^\ell {k_n \choose \ell} \Big(\frac{s}{n-1}\Big)^\ell \Big(1-\frac{s}{n-1}\Big)^{k_n -\ell}\Big)\notag \\
&= \Big(1-\frac{s}{n}\Big)^{k_n} \Big(1-\Big(1-\frac{s}{(n-s)(n-1)}\Big)^{k_n}\Big).
\end{align}
Using \eqref{ConditionalExpectation} and slightly adapting \eqref{PropCIEj},  the first part of $ T_2 $ can be handled in the following way:
\begin{align*}
&\EE\Big[\EE[\mathbf{1}_{E_{n,j}(k_n)\cap F_{n,j}^I}\vert \cC_n(k_n) ]\, \EE[\mathbf{1}_{E_{n,i}(k_n)\cap F_{n,i}^I}\vert \cC_n(k_n) ]\Big]\\
&=\EE\Big[ \mathbf{1}_{E_{n,i}(k_n)} \mathbf{1}_{E_{n,j}(k_n)} \Big(1-\big(1-\frac{1}{n-s}\big)^{\vert C_I\vert}\Big)^2\Big]\\
&= \sum_{\ell=0}^{k_n}  \Big(1-\big(1-\frac{1}{n-s}\big)^{\ell}\Big)^2\, \PP(\vert C_I\vert = \ell \vert E_{n,i}(k_n)\cap E_{n,j}(k_n))\, \PP(E_{n,i}(k_n)\cap E_{n,j}(k_n))\\
&=\Big(1-\frac{s}{n}\Big)^{k_n}\Big(1-\frac{s}{n-1}\Big)^{k_n} \sum_{\ell=0}^{k_n}  \Big(1-\big(1-\frac{1}{n-s}\big)^{\ell}\Big)^2\,{k_n \choose \ell} \Big(\frac{s}{n-2}\Big)^\ell \Big(1-\frac{s}{n-2}\Big)^{k_n -\ell}\\
&=\Big(1-\frac{s}{n}\Big)^{k_n}\Big(1-\frac{s}{n-1}\Big)^{k_n}\\
&\qquad\qquad\times\sum_{\ell=0}^{k_n}  \Big(1-2\big(1-\frac{1}{n-s}\big)^{\ell}+\big(1-\frac{1}{n-s}\big)^{2\ell}\Big)\,{k_n \choose \ell} \Big(\frac{s}{n-2}\Big)^\ell \Big(1-\frac{s}{n-2}\Big)^{k_n -\ell}\\
&= \Big(1-\frac{s}{n}\Big)^{k_n}\Big(1-\frac{s}{n-1}\Big)^{k_n}\\
&\qquad \qquad \times\Bigg(1-2\Big(1-\frac{s}{(n-2)(n-s)}\Big)^{k_n}+\Big(1-\frac{2s}{(n-2)(n-s)}+ \frac{s}{(n-2)(n-s)^2}\Big)^{k_n}\Bigg),
\end{align*}
independently of $ i $ and $ j $. Combining this with \eqref{ProbEjFj} yields
\begin{align*}
T_2&\leq n^2 \Big(1-\frac{s}{n}\Big)^{2k_n}\bigg[2 \bigg(\Big(1-\frac{s}{(n-s)(n-1)}\Big)^{k_n}-\Big(1-\frac{s}{(n-s)(n-2)}\Big)^{k_n} \bigg)\\
&\qquad \qquad \qquad \qquad + \Big(1- \frac{2s}{(n-s)(n-2)}+ \frac{s}{(n-2)(n-s)^2}\Big)^{k_n}- \Big(1-\frac{s}{(n-s)(n-1)}\Big)^{2k_n}\bigg].
\end{align*}
Now, there exists $ n_2=n_2(\alpha,s)  $ such that 
\begin{align*}
\Big(1-\frac{s}{(n-s)(n-1)}\Big)^{k_n}-\Big(1-\frac{s}{(n-s)(n-2)}\Big)^{k_n} \leq 2 \frac{\alpha \log(n)+x}{n^2}
\end{align*}
and
\begin{align*}
	\Big(1- \frac{2s}{(n-s)(n-2)}+ \frac{s}{(n-2)(n-s)^2}\Big)^{k_n}- \Big(1-\frac{s}{(n-s)(n-1)}\Big)^{2k_n} \leq 2 \frac{\alpha \log(n)+x}{n^2}
\end{align*}
holds for all $n\geq n_2$.
Combing this with \eqref{eq:ProbEmptyCellAlpha} we can conclude that there exists a constant $c_3(s,x) \in (0, \infty)$  such that for all $ n \geq \max\{3,n_2\} $,
\begin{align}\label{T2}
T_2 \leq c_3(s,x)\frac{\alpha \log(n)}{n^{2 \alpha}}.
\end{align}
Putting the bounds \eqref{T1} and \eqref{T2}  into \eqref{VarianceTerms} we see that there  exists a constant $ c_4(s,x) \in (0, \infty)$ such that
\begin{align}\label{ConditionalVariance}
\var\big(\EE[Z_{n,s}-Z^I_{n,s}\vert\cC_n(k_n) ]\big)\leq c_2(s,x) \frac{\alpha^2 \log(n)^2}{n^{1+\alpha}} + c_3(s,x)\frac{\alpha \log(n)}{n^{2 \alpha}} \leq c_4(s,x) \frac{\alpha \log(n)}{n^{2\alpha}},
\end{align}
for all $ n \geq \max\{3,n_2\} $, since $\alpha<1$ by assumption.

\begin{remark}\label{rm:OptimalVar}\rm 
Using \eqref{VarianceTerms} and the exact expressions for the probabilities appearing there, it can be shown that the order of \eqref{ConditionalVariance} is optimal in the sense that one can find another constant $\tilde{c}_4(s,x)\in(0,\infty)$ such that $\var(\EE[Z_{n,s}-Z^I_{n,s}\vert\cC_n(k_n) ])\geq \tilde{c}_4(s,x) \frac{\alpha \log(n)}{n^{2\alpha}}$ for sufficiently large $n$.
\end{remark}

\subsection{Bounding $ \EE[(Z_{n,s}-Z^I_{n,s})^2] $}
For the second term in \eqref{NormalBound} it remains to bound $ \EE[(Z_{n,s}-Z^I_{n,s})^2] $. We have
\begin{align}\label{SquaredDiff}
\EE[(Z_{n,s}-Z^I_{n,s})^2]&\leq \EE[\Big(\sum_{\substack{j=1\\ j \neq I}}^n \mathbf{1}_{E_{n,j}(k_n)\cap F_{n,j}^I}-1\Big)^2]\notag\\
&=\EE\Big[\Big(\sum_{\substack{j=1\\ j \neq I}}^n \mathbf{1}_{E_{n,j}(k_n)\cap F_{n,j}^I}\Big)^2\Big]-2\EE\Big[\sum_{\substack{j=1\\ j \neq I}}^n \mathbf{1}_{E_{n,j}(k_n)\cap F_{n,j}}\Big]+1\notag\\
&=   \EE\Big[\sum_{\substack{j=1\\ j \neq I}}^n \sum_{\substack{i=1 \\ i\neq j,I}}^n \mathbf{1}_{E_{n,j}(k_n)\cap F_{n,j}^I}  \mathbf{1}_{E_{n,i}(k_n)\cap F_{n,i}^I}\Big]-\EE\Big[\sum_{\substack{j=1\\ j \neq I}}^n \mathbf{1}_{E_{n,j}(k_n)\cap F_{n,j}^I}\Big]+1,
\end{align}
where the $-1$ in the first line comes from the artificially isolated cell with index $I$.
For the first sum note that
\begin{align*}
\PP(E_{n,i}(k_n)\cap F_{n,i}^I\cap E_{n,j}(k_n)\cap F_{n,j}^I)=\PP(F^I_{n,i}(k_n)\cap F_{n,j}^I\vert E_{n,j}(k_n)\cap E_{n,i}) \,\PP(E_{n,j}(k_n)\cap E_{n,j}),
\end{align*}
where $ \PP(E_{n,j}(k_n)\cap E_{n,j}(k_n)) $ can be bounded using \eqref{SBIndicatorProb} and \eqref{eq:ProbEmptyCellAlpha}.
To bound the remaining probability, we observe that similar to the considerations in \eqref{ProbEjFj} we have
\begin{align*}
&\PP( F_{n,j}^I\vert E_{n,i}(k_n) \cap E_{n,j}(k_n))\\
&= \sum_{\ell=0}^{k_n} \PP( F_{n,j}^I\vert E_{n,i}(k_n)\cap E_{n,j}(k_n)\cap\{ |C_I|=\ell \})\PP(|C_I|=\ell\vert E_{n,i}(k_n)\cap E_{n,j}(k_n) )\\
&= \sum_{\ell=0}^{k_n} {k_n \choose \ell} \Big(1-\big(1-\frac{1}{n-s}\big)^\ell\Big)\Big(\frac{s}{n-2}\Big)^\ell\Big(1-\frac{s}{n-2}\Big)^{k_n-\ell}\\
&=1-\Big(1-\frac{s}{(n-2)(n-s)}\Big)^{k_n},
\end{align*} 
independently of the choice of $i$ and $j$.
So we are left to deal with $ \PP( F_{n,i}^I\vert E_{n,i}(k_n)\cap E_{n,j}(k_n)\cap F_{n,j}^I) $. It holds that for $ i \neq j $
\begin{align*}
\PP( F_{n,i}^I\vert E_{n,j}(k_n) \cap E_{n,i}(k_n) \cap F_{n,j}^I)\leq   \PP( F_{n,i}^I\vert  E_{n,i}(k_n) \cap E_{n,j}(k_n)),
\end{align*}
since the event $ F_{n,j}^I $ implies that at least one particle from cell $ I $ is placed into cell $ C_j, $ reducing the chances of cell $ C_i $ receiving a particle. Hence,
\begin{align*}
\PP( F_{n,i}^I \cap F_{n,j}^I\vert E_{n,i}(k_n) \cap E_{n,j}(k_n)) \leq\Big(1-\Big(1-\frac{s}{(n-2)(n-s)}\Big)^{k_n}\Big)^2
\end{align*}
 and we finally obtain that for $ n \geq \max \{n_0,s\} $,
 \allowdisplaybreaks
\begin{align*}
&\EE\Big[\sum_{\substack{j=1\\ j\neq I}}^n \sum_{\substack{i=1 \\ i\neq j,I}}^n \mathbf{1}_{E_{n,j}(k_n)\cap F_{n,j}^I}  \mathbf{1}_{E_{n,i}(k_n)\cap F_{n,i}^I}\Big]\\
&= \sum_{\substack{j=1\\ j\neq I}}^n \sum_{\substack{i=1 \\ i\neq j, I}}^n \PP( F_{n,i}^I \cap F_{n,j}^I\vert E_{n,i}(k_n) \cap E_{n,j}(k_n))\PP( E_{n,i}(k_n) \cap E_{n,j}(k_n)) \\
&\leq (n-1)(n-2)  \Big(1-\Big(1-\frac{s}{(n-2)(n-s)}\Big)^{k_n}\Big)^2 \Big(1-\frac{s}{n}\Big)^{2k_n} \\
&\leq (n-1)(n-2) \Big(k_n \frac{s}{(n-2)(n-s)}\Big)^2 \Big(1-\frac{s}{n}\Big)^{2k_n}\\
&\leq (n-1)(n-2)\Big(\frac{3 (s+1)(\alpha \log(n)+x)}{n}\Big)^2 \Big( \frac{2 e^{-x}\alpha\log(n)}{n^{\alpha}}\Big)^2\\
&\leq 36(s+1)^2 e^{-2x}\frac{(\alpha\log(n)+x)^2}{n^{2 \alpha}},
\end{align*}
where we have used \eqref{Ineq:1-z^k} to arrive at the third line.
Since according to \eqref{ProbEjFj} the negative sum in \eqref{SquaredDiff} is of order $n^{-(\alpha+1)}$, we conclude that there exists a constant $ c_5(s,x)\in  (0, \infty) $  such that for all $  n \geq \max \{n_0,s\} $,
\begin{align}\label{ExpectationSquared}
	\EE[(Z_{n,s}-Z^I_{n,s})^2]\leq 
	c_5(s,x).
\end{align}
Combining the normal approximation bound \eqref{NormalBound} with the estimates \eqref{Expectation}, \eqref{Variance}, \eqref{ConditionalVariance} and \eqref{ExpectationSquared} we arrive at
\begin{align*}
d_W(\tilde{Z}_{n,s},G)&\leq \frac{\lambda_n}{c_1(s,x)\lambda_n}\sqrt{c_4(s,x) \frac{\alpha\log(n)}{n^{2\alpha}}}+\frac{\lambda_n}{(c_1(s,x)\lambda_n)^{\frac{3}{2}}}c_5(s,x)\\
&\leq  
\begin{cases}
C(s,x)\frac{\sqrt{\log(n)}}{n^{\alpha}} &\text{ for } \alpha \in (0, \frac{1}{3}]\\[0.2cm]
C(s,x)n^{-{1\over 2}(1-\alpha)}&\text{ for } \alpha \in( \frac{1}{3},1)
\end{cases}
\end{align*}
for some constant  $ C(s,x)\in (0, \infty ) $ and for all $n\geq N:=\max \{n_0,n_1,n_2,3,s, e^{-x}\}$. 
This proves Theorem~\ref{Thm:Normal}.\qed

\section{Proof of Theorem~\ref{Thm:Poisson}}\label{sec:PoissonProof}

As already mentioned in the introduction, rates of convergence towards a Poisson limit in the case $ \alpha=1$ can be concluded from \cite[Theorem~6.F]{BHS92}. Nevertheless, we give an alternative self-contained proof using the coupling from Section \ref{sec:Coupling}.

\medspace

It follows from \cite[Theorem 4.13]{Ross} that the total variation distance between the number of empty cells in the coupon collector's problem and a Poisson random variable $ W $ with  parameter $ \lambda_n=\EE[Z_{n,s}(k_n)] $ can be bounded using the following inequality:
\begin{align}\label{PoissonBound}
d_{TV}(Z_{n,s}(k_n), W) \leq \min\{1, \lambda_n\} \, \EE[Z_{n,s}(k_n)+1-Z_{n,s}^I(k_n)].
\end{align}
For the expectation on the right the definitions of $ Z_{n,s}(k_n) $ and $ Z^I_{n,s}(k_n) $ yield
\begin{align*}
\EE[Z_{n,s}(k_n)+1-Z_{n,s}^I(k_n)]&=\EE[\mathbf{1}_{E_{n,I}(k_n)}]+ \EE\Big[\sum_{j\neq I} \mathbf{1}_{E_{n,j}(k_n)}-\mathbf{1}_{E^I_{n,j}(k_n)}\Big]\\
&=\Big(1-\frac{s}{n}\Big)^k +(n-1)\bigg(\Big(1-\frac{s}{n}\Big)^k -\Big(1-\frac{s}{n-1}\Big)^k \bigg).
\end{align*}
 Combining \eqref{eq:ProbEmptyCellAlpha} and \eqref{eq:ProbEmptyCellsDiffAlpha} for $ \alpha =1 $  with  \eqref{PoissonBound} yields
\begin{align*}
d_{TV}(Z_{n,s}(k_n), W) &\leq \min\{1, \lambda_n\} \, \EE[Z_{n,s}(k_n)+1-Z_{n,s}^I(k_n)]\\
& = \Big(1-\frac{s}{n}\Big)^k +(n-1)\bigg(\Big(1-\frac{s}{n}\Big)^k -\Big(1-\frac{s}{n-1}\Big)^k \bigg)\\
&\leq \frac{2\, e^{-x}}{n}\Big(\log(n)+x+1\Big),
\end{align*}
for $ n \geq \tilde{N}:=\max\{ n_0, n_1, e^{-(x+1)}\} $,  which completes the proof of Theorem~\ref{Thm:Poisson}.\qed

\subsection*{Acknowledgement}

This work has been supported by the DFG priority program SPP 2265 \textit{Random Geometric Systems}.



\begin{thebibliography}{30}\small
	
\bibitem{ArratiaGoldsteinKochman}
R. Arratia, L. Goldstein, F. Kochman: Size bias for one and all. Probab. Surveys \textbf{16}, 1--61 (2019).
	
\bibitem{BHS92}
A.D. Barbour, L. Holst, S. Janson: \emph{Poisson Approximation}. Oxford University Press (1992).

\bibitem{BaumBillingsley}
L.E. Baum, P. Billingsley: Asymptotic distributions for the coupon collector’s problem. Ann. Math. Statistics \textbf{36}, 1835--1839 (1965).
	
\bibitem{Euler}
L. Euler: Solutio quarundam quaestionum difficiliorum in calculo probabilium. Opuscula Analytica, Vol. 2, 331--346 (1785).
	
	 
\bibitem{GR96}
L. Goldstein, Y. Rinott: Multivariate normal approximations by Stein's method and size bias couplings. J. Appl. Probab. \textbf{33}, 1--17 (1996).


\bibitem{Johnson}
N.L. Johnson, S. Kotz: \emph{Urn Models and their Application: An Approach to Modern Discrete Probability Theory}.  Wiley (1977).

\bibitem{Kolchin}
V.F. Kolchin,  B.A. Sevastianov, V.P. Chistiakov:   \textit{Random Allocations}.  Halsted Press Washington   (1978).


\bibitem{Laplace}
P.S. Laplace: \textit{Th\'eorie Analytique des Probabilit\'es}. Courcier, Paris (1812).


\bibitem{Mahmoud}
H. Mahmoud: Gaussian phases in generalized coupon collection. Adv. in Appl. Probab. \textbf{42}, 994-1012. (2010)

\bibitem{MikhailovPoisson77}
V.G. Mikhailov: A Poisson limit theorem in the scheme of group disposal of particels. Th. Probab. Appl. \textbf{22}, 152--156 (1977).

\bibitem{MikhailovNormal80}
V.G. Mikhailov: Asymptotic normality of the number of empty cells for group allocation of particles. Th. Probab. Appl. \textbf{25}, 82--90 (1980).

\bibitem{DeMoivre} 
A. de Moivre: \textit{The Doctrine of Chances}. W. Pearson (1718).

\bibitem{Ross} 
N. Ross: Fundamentals of Stein's method. Probab. Surveys \textbf{8}, 210--293 (2011).

\bibitem{SchillingHenze}
J. Schilling, N. Henze: Two Poisson limit theorems for the coupon collector's problem with group drawings. J. Appl. Probab. \textbf{58}, 966--977 (2021).

\bibitem{Smythe}
R. T. Smythe. Generalized coupon collection: the superlinear case. J. Appl. Probab. \textbf{48}, 189--199 (2011).

\bibitem{Stadje}
W. Stadje: The collector’s problem with group drawings. Adv. in Appl. Probab. \textbf{22}, 866--882 (1990).

\bibitem{VatutinMikhailov}
V.A. Vatutin, V.G. Mikhailov: Limit theorems for the number of empty cells in an equiprobable scheme for group allocation of particles. Th. Probab. Appl. \textbf{27}, 734--743 (1983).
		
\end{thebibliography}
\end{document}